\documentclass{amsart}
\usepackage{amssymb}

\newtheorem{theorem}{Theorem}

\theoremstyle{definition}

\theoremstyle{remark}

\numberwithin{equation}{section}

\begin{document}
\title{Linearization and local stability of random dynamical systems}
\author{Igor V. Evstigneev}
\address{Economics Department, University of Manchester, Manchester M13 9PL,
United Kingdom}
\email{igor.evstigneev@manchester.ac.uk}
\author{Sergey A. Pirogov}
\address{Institute for Information Transmission Problems, Academy of
Sciences of Russia, GSP-4, Moscow, 101447, Russia}
\email{pirogov@mail.ru}
\author{Klaus R. Schenk-Hopp\'{e}}
\address{School of Mathematics and Leeds University Business School,
University of Leeds, Leeds LS2 9JT, United Kingdom }
\email{k.r.schenk-hoppe@leeds.ac.uk}
\thanks{The authors gratefully acknowledge financial support from the Swiss
National Center of Competence in Research ``Financial Valuation and Risk
Management'' (project ``Behavioural and Evolutionary Finance'') and from the
Finance Market Fund, Norway (projects ``Stochastic Dynamics of Financial
Markets'' and ``Stability of Financial Markets: An Evolutionary Approach'').}
\date{March 29, 2010}

\begin{abstract}
The paper examines questions of local asymptotic stability of random
dynamical systems. Results concerning stochastic dynamics in general metric
spaces, as well as in Banach spaces, are obtained. The results pertaining to
Banach spaces are based on the linearization of the systems under study. The
general theory is motivated (and illustrated in this paper) by applications
in mathematical finance.
\end{abstract}

\maketitle

\section{Stochastic dynamics and equilibrium}\label{sec:1}

Let $(\Omega, \mathcal{F}, P)$ be a probability space and $T:\Omega
\rightarrow \Omega $ its endomorphism, i.e., an $\mathcal{F}$-measurable map
preserving the probability $P$:
\begin{equation*}
P(T^{-1}(\Gamma ))=P(\Gamma ),\;\Gamma \in \mathcal{F}.
\end{equation*}
Let $(X,\mathcal{X})$ be a measurable space and $f(x,\omega )$ a jointly
measurable mapping of $X\times \Omega $ into $X$. These data define a
discrete-time (autonomous) random dynamical system with the state space $X$,
the law of motion $f(x,\omega )$ and the time shift $T$. For each $t=1,2,...$, we put
\begin{equation}
f_{t}(x,\omega) := f(x,T^{t-1}\omega).  \label{f-t}
\end{equation}
A sequence $x_{0}(\omega), x_{1}(\omega), ...$ of random elements in $X$ is
called a \textit{path} of the random dynamical system if
\begin{equation}
x_{t}(\omega )=f_{t}(x_{t-1}(\omega ),\omega ), \; t=1,2,....  \label{l-m}
\end{equation}
almost surely (a.s.). A random element $\bar{x}(\omega )\in X$ satisfying
\begin{equation}
\bar{x}(T\omega ) = f(\bar{x}(\omega),\omega) \; \text{(a.s.)}  \label{equil}
\end{equation}
is termed a \textit{stochastic equilibrium} (or a \textit{random fixed point}, or a \textit{stationary point}) of the system. The sequence $\bar{x}_{t}(\omega) := \bar{x}(T^{t}\omega )$, which satisfies $\bar{x}_{t}(\omega) = f_{t}(\bar{x}_{t-1}(\omega ),\omega )$ (a.s.) by virtue of (\ref{equil}), represents the \textit{equilibrium path} of the system generated by the random fixed point $\bar{x}(\omega)$. For each $a \in X$, we denote by $x_{0}^{a}(\omega), x_{1}^{a}(\omega), ...$ the path with the initial state $a$ --- the sequence of random elements generated recursively (for each $\omega $ and all $t\geq 1$) by equations (\ref{l-m}) with $x_{0}(\omega) = a$.

\section{Basic assumptions}\label{sec:2}

Let $\bar{x}(\omega)$ be a stochastic equilibrium and let $%
X(\omega)\subseteq X$ be a random set\footnote{We say that $X(\omega)$ is a random set if the graph $\{(u,\omega)\in X \times \Omega: u\in X(\omega)\}$ of the multivalued mapping $\omega\mapsto X(\omega)$ belongs to the $\sigma $-algebra $\mathcal{X}\otimes\mathcal{F}$.} such that $\bar{x}(\omega)\in X(\omega)$ and $f(x,\omega)\in X(T\omega)$ for each $x\in X(\omega)$ and $\omega\in\Omega$. Let $\rho$ be a metric on $X$. For any number $\delta>0$,
define
\begin{equation*}
\bar{X}(\omega,\delta):=\{x\in X(\omega):\rho(x,\bar{x}(\omega))\leq\delta\}.
\end{equation*}
We introduce two hypotheses, (\textbf{A}) and (\textbf{B}), that will be
needed for the statement of the main results.

(\textbf{A}) $X$ is a complete separable metric space with the metric $\rho(x,y)$ and the Borel $\sigma$-algebra $\mathcal{X}$. There exist random
variables $L(\omega) > 0$ and $\delta(\omega) > 0$ such that
\begin{equation}
E|\ln L|<+\infty,\;E|\ln\delta|<+\infty,  \label{beta-delta}
\end{equation}
\begin{equation}
E\ln L<0,  \label{beta-0}
\end{equation}
and
\begin{equation}
\rho(f(x,\omega),f(\bar{x}(\omega),\omega))\leq L(\omega)\rho(x,\bar{x}%
(\omega))  \label{rho}
\end{equation}
for all $x\in\bar{X}(\omega,\delta(\omega))$. Here, ``$E$'' stands for the
expectation with respect to the probability measure $P$.

It is assumed in (\textbf{A}) that the mapping $f(\cdot,\omega)$ is locally
Lipschitzian at the equilibrium point $\bar{x}(\omega)$ (in a $\delta(\omega)$-neighborhood of $\bar{x}(\omega)$) with the random Lipschitz constant $L(\omega)$. Conditions (\ref{beta-delta}) guarantee that the random
Lipschitz constant is not ``too large'' and the random $\delta(\omega)$-neighborhood is not ``too small.'' (We could assume that $E \ln^{+} L < +\infty$ and $E\ln ^{-}\delta>-\infty$, but this would not lead to a gain in generality.) According to (\ref{beta-0}), $L(\omega)$ has negative expected logarithm, so that the mapping $f(\cdot,\omega)$ is \textit{locally contracting on average}.

The next hypothesis pertains to the case of a linear space $X$. It provides
conditions guaranteeing the validity of (\textbf{A}) formulated in terms of
the linearization of the mapping $f(\cdot,\omega)$ at the equilibrium point $\bar{x}(\omega)$.

(\textbf{B}) $X$ is a separable Banach space with the norm $\|\cdot\|$, the
metric $\rho(x,y) := \|x-y\|$ and the Borel $\sigma$-algebra $\mathcal{X}$.

(\textbf{B1) }There exist random variables $L(\omega) > 0$ and $\delta(\omega) > 0$ for which conditions (\ref{beta-delta}) hold and the mapping $f(x,\omega)$ satisfies (\ref{rho}) for all $x \in \bar{X}(\omega,\delta(\omega))$.

(\textbf{B2}) For each $\omega$, there exists the derivative $f^{\prime}(\bar{x}(\omega),\omega)$ of the mapping $f(\cdot,\omega)$ at the point $\bar{x}(\omega)$, and
\begin{equation*}
E\ln\|f^{\prime}(\bar{x}(\omega),\omega)\| < 0.
\end{equation*}

The derivative is understood in the strong (Fr\'{e}chet) sense, as a
continuous linear operator $F(\omega):=f^{\prime}(\bar{x}(\omega),\omega)$
such that
\begin{equation}
f(\bar{x}(\omega)+h,\omega) = f(\bar{x}(\omega),\omega) + F(\omega)h + g(h,\omega) \|h\|, \; h \in X,  \label{D}
\end{equation}
where $g(h,\omega) \rightarrow 0$ as $\|h\| \rightarrow 0$. The norm $\|F\|$ of the operator $F$ is defined as $\sup\{\|Fh\|/\|h\|: 0\neq h\in X\}$. Note
that hypotheses (\textbf{A}) and (\textbf{B1}), in contrast with (\textbf{B2}), impose assumptions on the behavior of the mapping $f(x,\omega)$ only in
the intersection of a neighborhood of the point $\bar{x}(\omega)$ with the
set $X(\omega)$ (which might be significantly smaller than the whole
neighborhood).

Throughout the paper we will assume that the endomorphism $T$ is ergodic:
all invariant sets have either zero or full measure.

\section{Central result}\label{sec:3}

\begin{theorem}
\label{thm:1} Let hypothesis (\textbf{A}) or hypothesis (\textbf{B}) hold.
Then there exist a random variable $\gamma(\omega) > 0$ and a constant $c < 0$ such that with probability one,
\begin{equation}
\lim\sup\nolimits_{t\rightarrow\infty}\frac{1}{t}\ln\rho(x_{t}^{a}(\omega),%
\bar{x}_{t}(\omega))\leq c  \label{gamma-l}
\end{equation}
for all $a\in\bar{X}(\omega,\gamma(\omega))$.
\end{theorem}

The theorem shows that for all $\omega$ in a set $\Omega_{1}$ of measure $1$, the distance between the path $x_{t}^{a}(\omega)$ with the initial state $a$ and the equilibrium path $\bar{x}_{t}(\omega)$ tends to zero at an
exponential rate for all $a$ in a sufficiently small neighborhood $\bar{X}(\omega,\gamma(\omega))$ of $\bar{x}(\omega)$. This convergence and its rate are uniform with respect to $\omega\in\Omega_{1}$ and $a\in\bar{X}(\omega,\gamma(\omega))$.

The results of this paper are motivated by applications in evolutionary
finance \cite{EvstigneevHensSchenkHoppe2009} --- a new and rapidly
developing area in financial mathematics. The general mathematical framework
for models in this area is the ergodic theory of random dynamical systems.
For the analysis of evolutionary models of asset markets, one often needs to
examine local stability of stochastic equilibria of some dynamical systems.
Surprisingly, the literature does not contain results which would be quite
suitable for applications in this field. The main tools offered for local
analysis in the context of stochastic dynamics (see the classical monograph
by L. Arnold \cite{Arnold1998}) are stochastic analogues of the
Hartman-Grobman theorem \cite{Grobman1959} -- \cite{Hartman1960a} developed
in \cite{Wanner1994} and \cite{CoaylaTeranMohammedRuffino2007} -- \cite{Cong1996}, and closely related results on stable manifold theorems \cite{Carverhill1985, LianLu2010, MohammedScheutzow1999, MohammedZhangZhao2008, Ruelle1979, Ruelle1982}. For the most part, these are delicate results, the
use of which requires the verification of complex conditions. They are much
stronger than what is needed for our purposes and require much stronger
assumptions. The need for suitable tools motivated us to develop the
techniques described above. To use Theorem~\ref{thm:1} under hypothesis (\textbf{B}) one has essentially to estimate only the expectation of one
random variable: the logarithm of the norm of the derivative $f^{\prime}(x,\omega)$ at the random fixed point $x = \bar{x}(\omega)$ (the local Lipschitz property (\textbf{B1}) usually causes no problems). The theorem represents a direct stochastic analogue of well-known deterministic results
on local asymptotic stability. Its statement and proof are based on
elementary notions and techniques.

By and large, stochastic equilibria, or random fixed points, exist under
sufficiently strong assumptions.\footnote{Counterexamples in \cite{EvstigneevPirogov2010, OchsOseledets1999} substantiate this assertion.} In many cases such assumptions guarantee not only existence, but also global asymptotic stability established simultaneously as a consequence of a stochastic contraction principle \cite{EvstigneevPirogov2007, Schmalfuss1998}. Global asymptotic stability --- convergence to a limit from any initial state --- implies that the limit forms an equilibrium. Of course in such cases, local analysis is superfluous. This is the situation, for example, with stochastic equilibrium problems related to random Perron-Frobenius theorems \cite{ArnoldGundlachDemetrius1994, Evstigneev1974, EvstigneevPirogov2010, Kifer1996b, Kifer2009}. Results in that field have been primarily inspired by the applications in the theory of Gibbs measures \cite{PirogovSinai19751976, Sinai1972}, large deviations \cite{Kifer1998}, models in mathematical biology \cite{DemetriusGundlach1999}, and others (see also \cite{KhaninKifer1996, Kifer1996a, KiferLyu2006}). In many models coming from economics and finance, the situation is different. The questions of existence are often separated from the questions of stability, and existence proofs based on the above methods require too restrictive, unjustifiable assumptions. Typically, either the system has an obvious, even deterministic, fixed point, or the existence is proved by methods quite distinct from the above-mentioned arguments based on global stability and limiting procedures (see, e.g., studies on von Neumann-Gale dynamics \cite{ArnoldEvstigneevGundlach1999, EvstigneevSchenkHoppe2008}). In the example we consider at the end of the paper, we deal with the former case: the system is random but the fixed point is deterministic.

The remainder of the paper is organized as follows. In section~\ref{sec:4}
we prove Theorem~\ref{thm:1}. At the end of the section, we provide a
condition sufficient for the validity of hypothesis (\textbf{B1}). Section~\ref{sec:5} gives an extension of Theorem~\ref{thm:1} applicable both to continuous and discrete time settings. Sections~\ref{sec:6} and \ref{sec:7}
analyze an example related to mathematical finance.

\section{Proof of Theorem~\protect\ref{thm:1}}\label{sec:4}

\textit{1st step.} Assume that hypothesis (\textbf{A}) holds. Replacing $\Omega$ by its subset of full measure, we can (in view of (\ref{equil}))
assume without loss of generality that the equations $\bar{x}_{t+1}(\omega) = f_{t}(\bar{x}_{t}(\omega),\omega)$ hold for all $t$ and all $\omega$,
rather than almost surely. Define $X_{t}(\omega ):=X(T^{t}\omega)$, $\delta_{t}(\omega) := \delta (T^{t}\omega )$, $L_{t}(\omega) := L(T^{t-1}\omega)$ and
\begin{equation*}
\bar{X}_{t}(\omega) := \{ x \in X_{t}(\omega) : \rho (x,\bar{x}_{t}(\omega)) \leq \delta_{t}(\omega)\}.
\end{equation*}
It follows from (\ref{rho}) that
\begin{equation*}
\sup_{x\in \bar{X}_{t-1}(\omega )}\rho (f_{t}(x,\omega),\bar{x}_{t}(\omega)) \leq L_{t}(\omega) \rho(x,\bar{x}_{t-1}(\omega)).
\end{equation*}
This implies the validity of the following assertion:

(*) If $x_{t-1}^{a}(\omega)\in\bar{X}_{t-1}(\omega)$, then
\begin{equation*}
\rho(x_{t}^{a}(\omega),\bar{x}_{t}(\omega))\leq
L_{t}(\omega)\rho(x_{t-1}^{a}(\omega),\bar{x}_{t-1}(\omega)).
\end{equation*}

We claim that there exists a random variable $\gamma(\omega) > 0$ for which
the inequalities
\begin{equation}
\gamma(\omega)\leq\delta_{0}(\omega),  \label{gamma}
\end{equation}
\begin{equation}
L_{t}(\omega)...L_{1}(\omega)\gamma(\omega)\leq\delta_{t}(\omega), \; t=1,2,...
\label{gamma1}
\end{equation}
hold with probability one. Indeed, these inequalities are satisfied if and
only if
\begin{equation}
\sigma(\omega):=\sup_{t=0,1,...}\frac{L_{t}(\omega) ... L_{0}(\omega)}{\delta_{t}(\omega)} \leq \frac{1}{\gamma(\omega)},  \label{sigma}
\end{equation}
where $L_{0}(\omega):=1$. It is sufficient to show that $P\{ \sigma < \infty
\} = 1$. Then relation (\ref{sigma}), and hence inequalities (\ref{gamma}) and
(\ref{gamma1}), hold a.s.\ for the random variable $\gamma$ defined as $\gamma(\omega)=1/\sigma(\omega)$ if $\sigma(\omega)<\infty$ and $\gamma(\omega)=1$ otherwise. To prove that $\sigma<\infty$ a.s., we use the
ergodic theorem and obtain that
\begin{equation*}
\frac{1}{t}\ln\frac{L_{t}...L_{0}}{\delta_{t}}=\frac{1}{t}\sum_{i=0}^{t}\ln
L_{i}-\frac{1}{t}\ln\delta_{t}\rightarrow E\ln L<0 \; \text{(a.s.)},
\end{equation*}
since $E|\ln L|<\infty$ and $E|\ln\delta|<\infty$. Consequently, $L_{t} ... L_{0} \delta_{t}^{-1}\rightarrow0$ a.s., which implies that $P\{\sigma < \infty\} = 1$.

Consider the random variable $\gamma(\omega) > 0$ constructed above and
denote by $\Omega_{1}$ the set of those $\omega$ for which inequalities (\ref{gamma}) and (\ref{gamma1}) hold and additionally,
\begin{equation}
\lim\nolimits_{t\rightarrow \infty} \frac{1}{t} \sum_{i=0}^{t} \ln L_{i} = E \ln L.  \label{beta3}
\end{equation}
This relation holds a.s.\ by virtue of the ergodic theorem, and so $P(\Omega
_{1})=1$. Take any $\omega \in \Omega _{1}$ and $a\in X(\omega )$ satisfying
$\rho (a,\bar{x}(\omega ))\leq \gamma (\omega )$. Let us show by induction
that for all $t\geq 0$,
\begin{equation}
x_{t}^{a}(\omega )\in \bar{X}_{t}(\omega),  \label{U-t}
\end{equation}
\begin{equation}
\rho (x_{t}^{a}(\omega),\bar{x}_{t}(\omega)) \leq L_{t}(\omega) ... L_{0}(\omega) \gamma(\omega).  \label{beta-a}
\end{equation}
For $t=0$, we have $a\in X_{0}(\omega) = X(\omega)$ by assumption and
\begin{equation}
\rho(a,\bar{x}(\omega)) \leq \gamma(\omega) \leq \delta_{0}(\omega)
\label{rho-delta}
\end{equation}
by virtue of (\ref{gamma}), so that $x_{0}^{a}(\omega) \in \bar{X}_{0}(\omega)$. Inequality (\ref{beta-a}) is true for $t=0$ in view of (\ref{rho-delta}) and because $\gamma(\omega) = L_{0}(\omega) \gamma(\omega)$. Suppose relations (\ref{U-t}) and (\ref{beta-a}) are valid for some $t$.
Then $x_{t+1}^{a}(\omega )\in X_{t+1}(\omega )$ because
\begin{equation*}
\rho(x_{t+1}^{a}(\omega),\bar{x}_{t+1}(\omega)) \leq L_{t+1}(\omega) \rho(x_{t}^{a}(\omega),\bar{x}_{t}(\omega)) \leq L_{t+1}(\omega) L_{t}(\omega) ... L_{0}(\omega) \gamma(\omega),
\end{equation*}
where the first inequality follows from assertion (*) and (\ref{U-t}), while
the second is a consequence of (\ref{beta-a}). By using (\ref{gamma1}) (with
$t+1$ in place of $t$), we conclude that $\rho(x_{t+1}^{a}(\omega), \bar{x}_{t+1}(\omega)) \leq \delta_{t+1}(\omega)$, which proves that the
analogues of relations (\ref{U-t}) and (\ref{beta-a}) hold for $t+1$.
Inequality (\ref{beta-a}) combined with (\ref{beta3}) implies (\ref{gamma-l}) with $c := E\ln L$.

\textit{2nd step.} To complete the proof of Theorem~\ref{thm:1} we show that
(\textbf{B}) implies (\textbf{A}). In this connection, we make some comments
regarding measurability. By applying (\ref{D}) with $k^{-1}h$ in place of $h$, letting $k \rightarrow \infty$, and using the joint measurability of $f(x,\omega)$, we obtain that the mappings $F(\omega )h$ and $g(h,\omega )$,
and hence the functions $\|F(\omega )h\|$ and $\|g(h,\omega )\|$, are
jointly measurable in $(h,\omega )$. The operator norm $F(\omega )$ depends
measurably on $\omega $ because $\|F(\omega )\|:=\sup_{k}\{\|F(\omega) h_{k}\| / \|h_{k}\|\}$, where $\{h_{k}\}$ is a countable dense subset in $X \backslash \{0\}$.

By using (\ref{D}) and (\ref{rho}), we obtain that for each $x \in \bar{X}(\omega,\delta(\omega))$,
\begin{equation*}
\|g(x-\bar{x}(\omega),\omega)\| \, \|x-\bar{x}(\omega)\| \leq \|f(x,\omega) - f(\bar{x}(\omega),\omega)\| + \| F(\omega)(x-\bar{x}(\omega)) \| \leq
\end{equation*}
\begin{equation*}
L(\omega) \|x-\bar{x}(\omega)\| + \|F(\omega)\| \, \|x-\bar{x}(\omega)\|,
\end{equation*}
which implies
\begin{equation}
\|g(x-\bar{x}(\omega),\omega)\|\leq L(\omega)+\|F(\omega)\|.  \label{g-L-F}
\end{equation}
Further, if $x\in\bar{X}(\omega,\delta(\omega))$, then
\begin{equation}
\|f(x,\omega)-f(\bar{x}(\omega),\omega)\| \leq [ \|F(\omega)\| + \| g(x-\bar {x}(\omega),\omega)\| ] \, \|x-\bar{x}(\omega)\|.  \label{f-x}
\end{equation}
Define
\begin{equation}
g_{k}(\omega) := \sup \{ \|g(x-\bar{x}(\omega),\omega)\| : x \in \bar{X}(\omega,k^{-1}\delta(\omega)) \}.  \label{psi}
\end{equation}
The function $g_{k}(\omega)$ is measurable with respect to the completion $\mathcal{F}^{P}$ of the $\sigma$-algebra $\mathcal{F}$ by $P$-null sets
because $\|g(x-\bar{x}(\omega),\omega)\|$ is jointly measurable in $(x,\omega)$ and $X(\omega)$ is a random set. This follows from the fact that
the projection of a set in $\mathcal{X} \otimes \mathcal{F}$ on $\Omega$ is $\mathcal{F}^{P}$-measurable (see, e.g., \cite{DellacherieMeyer1978}, Theorem
III.33).

Define $D_{k}(\omega) := \|F(\omega )\| + g_{k}(\omega)$. By virtue of (\ref{D}), $g_{k}(\omega) \rightarrow 0$ for each $\omega$. Furthermore, in view of (\ref{g-L-F}) we have
\begin{equation*}
\ln D_{k}(\omega) \leq \ln [ 2 \|F(\omega)\| + L(\omega) ] \leq \ln 4 + \max [ \ln \|F(\omega )\|, \ln L(\omega)] =: \Xi(\omega),
\end{equation*}
where $E|\Xi (\omega)| < \infty$. By using Fatou's lemma, we get
\begin{equation*}
\lim \sup\nolimits_{k \rightarrow \infty} E\ln D_{k}(\omega) \leq
E\lim\nolimits_{k \rightarrow \infty} \ln D_{k}(\omega) = E \ln \|F(\omega)\| < 0.
\end{equation*}
Thus there exists $k$ such that $E\ln D_{k}(\omega )<0$ and (by virtue of (\ref{f-x}) and (\ref{psi}))
\begin{equation*}
\|f(x,\omega) - f(\bar{x}(\omega),\omega)\| \leq D_{k}(\omega) \|x - \bar{x}(\omega)\|
\end{equation*}
for each $x\in \bar{X}(\omega,\delta(\omega)/k)$. This completes the
proof of Theorem~\ref{thm:1}.\hfill $\Box$
\medskip

\noindent \textbf{Remark 1.} The following assumption is sufficient for
condition (\textbf{B1}) to hold.

(\textbf{B3}) There exist random variables $L(\omega ) > 0$ and $\delta(\omega) > 0$ satisfying (\ref{beta-delta}) such that for all $x \in \bar{X}(\omega,\delta(\omega))$, the set $X(\omega)$ contains the segment $[\bar{x}(\omega),x]$ connecting $\bar{x}(\omega)$ and $x$, the mapping $f(x,\omega)$ is differentiable at the point $x$, and the norm of the
derivative $\|f^{\prime }(x,\omega)\|$ is bounded from above by $L(\omega)$.

To deduce (\textbf{B1}) from (\textbf{B3}) it suffices to observe that the
inequality $\|f^{\prime}(x,\omega)\| \leq L(\omega)$ implies (\ref{rho}) by
virtue of the generalized mean value theorem:
\begin{equation*}
\|f(x,\omega) - f(\bar{x}(\omega),\omega)\| \leq \|x-\bar{x}(\omega)\| \sup_{y \in [\bar{x}(\omega),x]} \|f^{\prime}(y,\omega)\|
\end{equation*}
(which holds even if $f^{\prime}(x,\omega)$ is the weak rather than the
strong derivative) --- see \cite{KolmogorovFomin1957}, Section X.1.3.

\section{Extension to discrete and continuous time cocycles}\label{sec:5}

The next result provides an extension of Theorem~\ref{thm:1} to random
dynamical systems defined in terms of cocycles in discrete and continuous
time. Let $\mathbb{T}$ be either the set of non-negative integers or the set
of non-negative real numbers, and let $T^{t},t\in\mathbb{T},$ be a semigroup
of ergodic endomorphisms of the probability space $(\Omega,\mathcal{F},P)$.
For each $t\in\mathbb{T}$ let $C_{t}(x,\omega)$ be a jointly measurable
mapping of $X\times\Omega$ into $X$. Assume that the family of mappings $C_{t}(x,\omega)$, $t\in\mathbb{T}$, forms a \textit{cocycle}, i.e.,
\begin{equation*}
C_{t+s}(\cdot,\omega) = C_{s}(\cdot,T^{t}\omega) \circ C_{t}(\cdot,\omega), \; C_{0}(x,\omega)=x,
\end{equation*}
for all $t,s,x$ and $\omega$. The cocycle defines the law of motion in the
system, whose \textit{paths} are random functions $x_{t}(\omega)$, $t\in
\mathbb{T}$, such that with probability one, $x_{t}(\omega) = C_{t}(x_{0}(\omega),\omega)$ for all $t\in\mathbb{T}$. Random dynamical systems of this kind can be generated by stochastic or random differential equations in continuous time and by products of random mappings in discrete time (see \cite{Arnold1998}). We will assume that a random set $X(\omega)$ is given such that $C_{t}(x,\omega)\in X(T^{t}\omega)$ for all $x\in X(\omega)$.

Let $\bar{x}_{t}(\omega)$, $t \in \mathbb{T}$, be an \textit{equilibrium path}, i.e., a path satisfying a.s.\ $\bar{x}_{t}(\omega)=\bar{x}_{0}(T^{t}\omega)$
for all $t\in\mathbb{T}$. For each $a\in X$, define $x_{t}^{a}(\omega) := C_{t}(a,\omega)$ (the random path with the initial state $a$).

\begin{theorem}\label{thm:2}
Let the following assumptions hold:

(\textbf{C1}) There exists $M \in \mathbb{T}$ such that the mapping $C_{M}(\cdot,\omega)$ satisfies condition (\textbf{A}) or (\textbf{B}) with $\bar{x}(\omega) = \bar{x}_{0}(\omega)$.

(\textbf{C2}) There are random variables $H(\omega) > 0$, $\kappa(\omega) >
0 $ and a constant $b > 0$ such that $E|\ln H(\omega)| < +\infty$, $E|\ln\kappa(\omega)| < +\infty$ and with probability one,
\begin{equation}
\rho(C_{t}(x,\omega),C_{t}(\bar{x}(\omega),\omega)) \leq H(\omega) \rho(x,\bar{x}(\omega))^{b} \label{Holder}
\end{equation}
for all $x\in\bar{X}(\omega,\kappa(\omega))$ and $t \in \mathbb{T}$ satisfying
$0\leq t\leq M$.

Then there exist a random variable $\gamma(\omega) > 0$ and a constant $d <
0$ such that almost surely
\begin{equation}
\lim\sup\nolimits_{\mathbb{T} \ni t \rightarrow \infty}\frac{1}{t}\ln\rho
(x_{t}^{a}(\omega),\bar{x}_{t}(\omega))\leq d \label{gamma-m}
\end{equation}
for all $a \in \bar{X}(\omega,\gamma(\omega))$.
\end{theorem}

Since $d < 0$, (\ref{gamma-m}) implies that $\rho(x_{t}^{a}(\omega), \bar{x}_{t}(\omega)) \rightarrow 0$ a.s.\ at an exponential rate. This convergence, as well as its rate, are uniform with respect to $\omega$ in a set $\Omega_{1}$ of full measure and all $a$ in the neighborhood $\bar{X}(\omega,\gamma(\omega))$ of $\bar{x}(\omega) = \bar{x}_{0}(\omega)$. Property (\ref{Holder}) represents a H\"{o}lder condition on the cocycle $C_{t}(x,\omega)$.

\textit{Proof of Theorem~\ref{thm:2}.} Since (\textbf{B}) implies (\textbf{A}), as we demonstrated in section~\ref{sec:4}, it is sufficient to prove the
theorem under assumption (\textbf{A}). Define $t(n) := nM$, $n=0,1,...$. By
applying Theorem~\ref{thm:1} to the mapping $f(x,\omega) := C_{M}(x,\omega)$, we obtain that there exist a random variable $\gamma(\omega) > 0$ and a
constant $c<0$ such that with probability one,
\begin{equation}
\lim \sup\nolimits_{n\rightarrow \infty} \frac{1}{t(n)} \ln \rho(x_{t(n)}^{a}(\omega), \bar{x}_{t(n)}(\omega)) \leq c \label{gamma-l-M}
\end{equation}
for all $a \in \bar{X}(\omega,\gamma(\omega))$. From (\ref{gamma-l-M}) we
obtain that for each $\varepsilon \in (0,1)$ with probability one,
\begin{equation}
\rho _{n}(\omega ):=\rho(x_{t(n)}^{a}(\omega), \bar{x}_{t(n)}(\omega)) \leq
\exp [t(n)c(1-\varepsilon)]  \label{ineq1}
\end{equation}
for all sufficiently large $n$. Furthermore, almost surely
\begin{equation}
\kappa_{n}(\omega) := \kappa(T^{t(n)}\omega) \geq \exp [t(n) c (1-\varepsilon)]  \label{ineq2}
\end{equation}
for all $n$ large enough. Indeed, the last inequality holds if $M c (1-\varepsilon) \leq n^{-1} \ln \kappa_{n}(\omega)$, which is true for
all $n$ large enough because $\lim\nolimits_{n \rightarrow \infty} n^{-1} \ln
\kappa_{n}(\omega) \rightarrow 0$ a.s.\ (this follows from the assumption $E|\ln \kappa(\omega)| < +\infty $). From (\ref{ineq1}) and (\ref{ineq2}) we
obtain that the inequalities
\begin{equation}
\rho_{n}(\omega) \leq \kappa_{n}(\omega)  \label{rho-kappa}
\end{equation}
a.s.\ hold for all $n$ large enough. By using (\textbf{C2}) and (\ref{rho-kappa}), we obtain that for each $\varepsilon > 0$ with probability one,
\begin{equation*}
\eta _{n}(\omega ):=\sup_{t(n)\leq t\leq t(n+1)}\frac{1}{t}\ln \rho
(x_{t}^{a}(\omega ),\bar{x}_{t}(\omega ))\leq
\end{equation*}
\begin{equation*}
\sup_{t(n)\leq t \leq t(n+1)} \Big[ \frac{\ln^{+} H(T^{t(n)}\omega)}{t} + \frac{b}{t} \ln \rho(x_{t(n)}^{a}(\omega), \bar{x}_{t(n)}(\omega)) \Big] \leq
\end{equation*}
\begin{equation*}
\frac{\ln^{+} H(T^{t(n)}\omega)}{t(n)}+\sup_{t(n) \leq t \leq t(n+1)} \frac{b t(n) c (1-\varepsilon)}{t} = \frac{\ln^{+} H(T^{t(n)}\omega)}{t(n)} + \frac{b n c (1-\varepsilon)}{n+1}
\end{equation*}
for all $n$ greater than some $n(\omega)$. The first summand in the last
expression converges to zero a.s.\ by virtue of the assumption $E|\ln
H| < +\infty$. The second summand tends to $b c (1-\varepsilon)$.
Consequently, with probability one, there exists $k(\omega)$ such that $\eta_{n}(\omega) \leq b c (1-\varepsilon)^{2}$ for all $n\geq k(\omega)$.
Denote by $N(t)$ the natural number such that $M N(t) \leq t < M[N(t)+1]$. Then $t^{-1}\ln \rho (x_{t}^{a}(\omega),\bar{x}_{t}(\omega)) \leq \eta_{N(t)}(\omega)$. Thus if $t \geq M k(\omega)$, then $N(t)\geq k(\omega)$, and so
\begin{equation*}
t^{-1} \ln \rho(x_{t}^{a}(\omega), \bar{x}_{t}(\omega)) \leq \eta_{N(t)}(\omega) \leq b c (1-\varepsilon )^{2}.
\end{equation*}
Since $\varepsilon$ is any number in $(0,1)$, we obtain that (\ref{gamma-m}) holds with $d:=bc$ ($<0$). The proof is complete.\hfill $\Box$\medskip

\noindent\textbf{Remark 2.} Applying Theorem~\ref{thm:2} in the
discrete-time case, where $\mathbb{T} = \{0,1,...\}$, to the cocycle defined
by
\begin{equation*}
C_{t}(\cdot, \omega) := f_{t}(\cdot, \omega) \circ ... \circ f_{1}(\cdot, \omega), \; t\geq 1,
\end{equation*}
we obtain a version of Theorem~\ref{thm:1} in which condition (\textbf{A})
or (\textbf{B}) is imposed not on the given mapping $f(\cdot ,\omega)$, but
on the product $C_{M}(\cdot,\omega) = f_{M}(\cdot,\omega) \circ ...\circ
f_{1}(\cdot,\omega)$ of the mappings $f_{t}(\cdot,\omega)$. In this
case, hypothesis (\textbf{C2}) is fulfilled under the following assumption.

(\textbf{C3}) There exist random variables $L(\omega )>0$ and $\delta
(\omega )>0$ for which conditions (\ref{beta-delta}) and (\ref{rho}) hold.

Indeed, assume that (\textbf{C3}) is satisfied and put
\begin{equation*}
\kappa(\omega) := \min_{0 \leq t \leq M} \frac{\delta_{t}(\omega)}{L_{0}(\omega) ... L_{t}(\omega)},
\end{equation*}
where $L_{t}$ and $\delta_{t}$ are defined in section~\ref{sec:4}. Arguing
by induction and using (\ref{rho}), we obtain that
\begin{equation*}
\rho(C_{t}(x,\omega), C_{t}(\bar{x}(\omega),\omega)) \leq L_{0}(\omega) ... L_{t}(\omega) \rho(x,\bar{x}(\omega)) \leq \delta_{t}(\omega)
\end{equation*}
($t=1,2,...,M$), as long as $\rho(x,\bar{x}(\omega)) \leq \kappa(\omega)$. This yields (\textbf{C2}) with $b=1$ and $H(\omega) := L_{1}(\omega) ... L_{M}(\omega)$.\medskip

\noindent \textbf{Remark 3.} Assume that $X$ is a separable Banach space and
for each $\omega $, the mapping $C_{t}(\cdot ,\omega )$ is differentiable at
the point $\bar{x}_{0}(\omega )$ and continuous in a neighborhood of this
point. Then the family of linear operators $F_{t}(\omega )=C_{t}^{\prime }(%
\bar{x}_{0}(\omega ),\omega )$ forms a cocycle (this follows from the chain
rule of differentiation) --- the linearization of the cocycle $C_{t}(x,\omega )$. Suppose that $E\ln ^{+}\|F_{t}(\omega )\|<+\infty $. In
this case, the assumption that $E\ln \|F_{M}(\omega )\|<0$ for some $M>0$,
needed for the application of Theorem~\ref{thm:2} under hypothesis (\textbf{B}), is equivalent to the assumption that the Furstenberg-Kesten constant \cite{FurstenbergKesten1960, Oseledets1968, Arnold1998}
\begin{equation*}
\lim\nolimits_{t\rightarrow \infty }\frac{1}{t}E\ln \|F_{t}(\omega)\| = \inf_{t>0}\frac{1}{t}E\ln \|F_{t}(\omega )\|
\end{equation*}
is negative.

\section{Application to an investment model}\label{sec:6}

To describe an example to which we will apply Theorem~\ref{thm:1}, assume
that together with the probability space $(\Omega,\mathcal{F},P)$ and its
endomorphism $T$, we are given a family of $\sigma$-algebras $\mathcal{F}%
_{0}\subseteq\mathcal{F}_{1}\subseteq...\subseteq\mathcal{F}$ such that $%
T^{-1}(\Gamma)\in\mathcal{F}_{t+1}$ if and only if $\Gamma\in\mathcal{F}_{t}$
($\mathcal{F}_{t}$ contains events observable prior to time $t$). Denote by $%
\Delta$ the unit simplex $\{v=(v_{1},...,v_{K})\geq0:\sum v_{k}=1\}$. Let $%
R_{t}(\omega)=R(T^{t}\omega)$, $\lambda_{t}(\omega)=\lambda(T^{t}\omega)$
and $\lambda_{t}^{\ast}(\omega)=\lambda^{\ast}(T^{t}\omega)$ be stationary
processes with values in $\Delta$ adapted to the filtration $(\mathcal{F}%
_{t})$ and $r$ a number in $(0,1)$ such that
\begin{equation}
rE_{t}\lambda^{\ast}(T\omega)+(1-r)E_{t}R(T\omega)=\lambda^{\ast}(\omega)\;
\text{(a.s.)},  \label{G-Kelly}
\end{equation}
where $E_{t}(\cdot)=E(\cdot|\mathcal{F}_{t})$. The existence and uniqueness
of the solution $\lambda^{\ast}(\cdot)$ to equation (\ref{G-Kelly}) follows
from the Banach contraction principle and the fact that $r<1$.

Consider the random dynamical system whose paths $(x_{t})$ are defined by
\begin{equation}
x_{t+1}=x_{t}\frac{\sum_{k=1}^{K}[r\lambda_{t+1,k}^{\ast}+(1-r)R_{t+1,k}]%
\dfrac{\lambda_{t,k}}{\lambda_{t,k}x_{t}+\lambda_{t,k}^{\ast}}}{\sum
_{k=1}^{K}[r\lambda_{t+1,k}+(1-r)R_{t+1,k}]\dfrac{\lambda_{t,k}^{\ast}}{%
\lambda_{t,k}x_{t}+\lambda_{t,k}^{\ast}}},  \label{f-fin}
\end{equation}
where $R_{t,k}$, $\lambda_{t,k}$, and $\lambda_{t,k}^{\ast}$ are the
coordinates of the vectors $R_{t}$, $\lambda_{t}$, and $\lambda_{t}^{\ast}$, respectively. In the evolutionary model of an asset market developed in
\cite{EvstigneevHensSchenkHoppe2006, EvstigneevHensSchenkHoppe2009} (see
these references for details), $R_{t,k}\,$are \textit{relative dividends} of
$K$ assets, the sequences of vectors $\lambda=(\lambda_{t})$ and $\lambda^{\ast} = (\lambda_{t}^{\ast})$ are \textit{investment strategies} (\textit{portfolio rules}) and $r$ is the \textit{investment rate}. The
vectors $\lambda_{t}\in\Delta$ and $\lambda_{t}^{\ast}\in\Delta$ indicate
proportions according to which investors using the strategies $\lambda$ and $\lambda^{\ast}$ allocate wealth across assets. The strategy $\lambda^{\ast}$
defined by (\ref{G-Kelly}) is a generalization of the \textit{Kelly
portfolio rule}, well-known in mathematical finance (see, e.g., \cite{MacLeanThorpZiemba2011}). It is assumed that there are two groups of
investors, one using the strategy $\lambda^{\ast}$ and the other any
strategy $\lambda$ distinct from $\lambda^{\ast}$. The variable $x_{t}$
represents the ratio $w_{t}/w_{t}^{\ast}$ where $w_{t}^{\ast}$ and $w_{t}$
denote the total wealth of the former and the latter groups of investors,
respectively. The local stability of the dynamical system under
consideration at the fixed point $\bar{x}=0$ means that the portfolio rule $%
\lambda^{\ast}$ is \textit{evolutionary stable}. If the initial relative
wealth $x_{0}=a$ of the $\lambda$-investors (``mutants''--- in the
terminology borrowed from biology) is small enough, then they will be
eventually driven out of the market by the $\lambda^{\ast}$-investors: their
relative wealth $x_{t}^{a}$ will tend to zero.

Define $\lambda_{k}^{\ast} := \lambda_{0,k}^{\ast}$ and assume that the
following conditions hold.

(\textbf{K1}) $E\ln\min_{k}\lambda_{k}^{\ast} > -\infty$.

(\textbf{K2}) The random variables $\mu_{k}:=r\lambda_{1,k}^{\ast
}+(1-r)R_{1,k}$, $k=1,...,K$, are conditionally linearly independent given $\mathcal{F}_{0}$, i.e., the equality $\alpha_{1} \mu_{1} + ... + \alpha_{K} \mu_{K} = 0$ holding (a.s.) for some $\mathcal{F}_{0}$-measurable random variables $\alpha_{k}$ implies $\alpha_{1} = ... = \alpha_{K}=0$ (a.s.).

\begin{theorem}
\label{thm:3} For any strategy $\lambda \neq \lambda^{\ast}$, there exist a
random variable $\gamma(\omega) > 0$ and a constant $c < 0$ such that with
probability one, $\lim\sup\nolimits_{t\rightarrow\infty} t^{-1}\ln
x_{t}^{a}\leq c$ for all $0\leq a\leq\gamma(\omega)$.
\end{theorem}

\textit{Proof.} We apply Theorem~\ref{thm:1} with $X=(-\infty,+\infty)$ and $X(\omega)=[0,\infty)$. We define the function $f(x,\omega)$ as the
right-hand side of (\ref{f-fin}) with $t=0$ if $x=x_{t}>-\zeta(\omega)$,
where $\zeta(\omega):=\min_{k}\lambda_{k}^{\ast}(\omega)$. For $x \leq -\zeta(\omega )$, we can define $f(x,\omega)$, for example, as any
(fixed) number $u$. Clearly $\bar{x}:=0$ is a fixed point of $f(x,\omega)$
for each $\omega$. For $0<x\leq1$, we have $0<x^{-1}f(x,\omega)\leq 2 \zeta^{-2}$, so that (\textbf{B1}) holds with $\delta=1$ and $L=2\zeta^{-2}$. Further,
\begin{equation*}
f^{\prime}(0,\omega)=\sum_{k=1}^{K}\mu_{k}\dfrac{\lambda_{k}}{\lambda
_{k}^{\ast}}.
\end{equation*}
By virtue of Jensen's inequality, we have
\begin{equation*}
E\ln\sum_{k=1}^{K}\mu_{k}\dfrac{\lambda_{k}}{\lambda_{k}^{\ast}}%
=E(E_{0}\ln\sum_{k=1}^{K}\mu_{k}\dfrac{\lambda_{k}}{\lambda_{k}^{\ast}}%
)<E(\ln \sum_{k=1}^{K}E_{0}\mu_{k}\dfrac{\lambda_{k}}{\lambda_{k}^{\ast}}%
)=E\ln \sum_{k=1}^{K}\lambda_{k}=0.
\end{equation*}
To show that the above inequality is indeed strict, assume the contrary.
Then the random variable $\sum_{k=1}^{K}\mu_{k}\lambda_{k}/\lambda
_{k}^{\ast}$ coincides a.s.\ with an $\mathcal{F}_{0}$-measurable random
variable. Hence it coincides a.s.\ with its conditional expectation given $\mathcal{F}_{0}$, which is equal to $1$ (this follows from (\ref{G-Kelly})).
Thus $\sum_{k=1}^{K}\mu_{k} \lambda_{k}/\lambda_{k}^{\ast} = 1$ (a.s.) or
equivalently, $\sum_{k=1}^{K}\mu_{k}[(\lambda_{k}/\lambda_{k}^{\ast })-1] = 0$
(a.s.), which implies by virtue of (\textbf{K2}) that $(\lambda_{k}/\lambda_{k}^{\ast}) - 1 = 0$ (a.s.) for all $k=1,2,...,K$. Consequently, $\lambda=\lambda^{\ast}$, which is a contradiction.\hfill $\Box$

\section{Sufficient conditions in the Markovian case}\label{sec:7}

We conclude with some comments regarding hypotheses (\textbf{K1}) and (\textbf{K2}). These hypotheses are formulated in terms of the Kelly strategy
$\lambda^{\ast}$, which is defined as the solution to equation (\ref{G-Kelly}). For the applications it is important to provide conditions sufficient for
(\textbf{K1}) and (\textbf{K2}) that are formulated in terms of one of the
primitives of the model --- the dividend process $R_{t}$. The former
condition holds if $E \ln \min_{k} E_{0} R_{1,k} > -\infty$, which is clear from (\ref{G-Kelly}). The latter is satisfied, for example, if the following
requirements are fulfilled:

(i) the random variables $R_{1,k}$, $k=1,...,K$, are conditionally
independent given $\mathcal{F}_{0}$ (the absence of ``redundant'' assets);

(ii) the probability space $(\Omega,\mathcal{F},P)$ and the filtration $(\mathcal{F}_{t})$ are generated by a stationary Markov process $...,s_{-1},s_{0},s_{1},...$ with values in a measurable space $S$, and the
vector function $R(\omega)$ depends only on $s_{0}$ (we shall denote it as $R(s_{0})$);

(iii) the process $s_{t}$ has a transition function $\pi(s,d\sigma)$
possessing a jointly measurable density $p(s,\sigma)$ with respect to a
probability measure $\pi(d\sigma)$ such that $0<v\leq p(s,\sigma)\leq V$ for
some constants $v \leq V$.

Let us prove that conditions (i)--(iii) imply (\textbf{K2}). We first
observe that in the Markov case, $\lambda_{k}^{\ast} = \lambda_{k}^{\ast}(s_{0})$ and $\mu_{k}=\mu_{k}(s_{1})$ are functions of $s_{0}$ and $s_{1}$, respectively. It follows from (\ref{G-Kelly}) that $\mu$ satisfies
\begin{equation}
rE_{0}\mu(s_{1})+(1-r)R(s_{0})=\mu(s_{0}) \; \text{(a.s.)}.  \label{mu-R}
\end{equation}
In the present setting, $\mathcal{F}_{t}$-measurable functions can be
identified with measurable functions $\alpha(s^{t})$, where $s^{t} := (..., s_{t-1}, s_{t})$. Let $\alpha_{1}(s^{0}), ..., \alpha_{K}(s^{0})$
be an $\mathcal{F}_{0}$-measurable vector functions satisfying
\begin{equation}
\langle\alpha(s^{0}),\mu(s_{1})\rangle=0 \; \text{(a.s.)},  \label{gamma-mu}
\end{equation}
where $\alpha:=(\alpha_{1},...,\alpha_{K})$ and $\mu=(\mu_{1},...,\mu_{K})$. We have to prove that $\alpha=0$ (a.s.). From (\ref{mu-R}) we get
\begin{equation}
r\langle\alpha(s^{-1}),E_{0}\mu(s_{1})\rangle+(1-r)\langle\alpha
(s^{-1}),R(s_{0})\rangle=\langle\alpha(s^{-1}),\mu(s_{0})\rangle \; \text{(a.s.)}.  \label{mu-R1}
\end{equation}
Let us show that $E|\langle\alpha(s^{-1}),E_{0}\mu(s_{1})\rangle|=0$. We
have
\begin{equation*}
E |\langle\alpha(s^{-1}), E_{0}\mu(s_{1}) \rangle| = E | E_{0} \langle \alpha(s^{-1}), \mu(s_{1}) \rangle | \leq E E_{0} | \langle \alpha(s^{-1}), \mu(s_{1})\rangle | =
\end{equation*}
\begin{equation}
E | \langle \alpha(s^{-1}), \mu(s_{1}) \rangle | = E E_{-1} |\langle \alpha(s^{-1}), \mu(s_{1}) \rangle|.  \label{mu-R5}
\end{equation}
From (iii) we get $p(s,\sigma)\leq Vv^{-1}p(s_{-1},\sigma)$. By using this,
we obtain
\begin{equation*}
E_{-1}|\langle\alpha(s^{-1}),\mu(s_{1})\rangle|=\int\pi(s_{-1},ds)\int
\pi(s,d\sigma)|\langle\alpha(s^{-1}),\mu(\sigma)\rangle|=
\end{equation*}
\begin{equation*}
\int\pi(s_{-1},ds)\int
p(s,\sigma)\pi(d\sigma)|\langle\alpha(s^{-1}),\mu(\sigma)\rangle|\leq
\end{equation*}
\begin{equation*}
Vv^{-1}\int\pi(s_{-1},ds)\int p(s_{-1},\sigma)\pi(d\sigma)|\langle
\alpha(s^{-1}),\mu(\sigma)\rangle|=
\end{equation*}
\begin{equation}
Vv^{-1}\int\pi(s_{-1},d\sigma)|\langle\alpha(s^{-1}),\mu(\sigma)\rangle
|=Vv^{-1}E_{-1}|\langle\alpha(s^{-1}),\mu(s_{0})\rangle|.  \label{mu-R7}
\end{equation}
By combining (\ref{mu-R7}), (\ref{mu-R5}) and (\ref{gamma-mu}), we obtain
that $E|\langle\alpha(s^{-1}),E_{0}\mu(s_{1})\rangle|=0$. This, together
with (\ref{gamma-mu}) and (\ref{mu-R1}), implies the equality $\langle \alpha(s^{-1}), R(s_{0}) \rangle = 0$ (a.s.). By using (i) we conclude that $\alpha(s^{-1})=0$ (a.s.), or equivalently, $\alpha(s^{0})=0$ (a.s.), which
completes the proof.


\begin{thebibliography}{99}
\bibitem{Arnold1998} L. Arnold. \textit{Random Dynamical Systems}. Springer,
1998.

\bibitem{ArnoldEvstigneevGundlach1999} L. Arnold, I. V. Evstigneev and V. M.
Gundlach. Convex-valued random dynamical systems: A variational principle
for equilibrium states. \textit{Random Oper. Stoch. Equ.} \textbf{7}
(1999), 23--38.

\bibitem{ArnoldGundlachDemetrius1994} L. Arnold, V. M. Gundlach and L.
Demetrius. Evolutionary formalism for products of positive random matrices.
\textit{Ann. Appl. Probab.} \textbf{4} (1994), 859--901.

\bibitem{Carverhill1985} A. Carverhill. Flows of stochastic dynamical
systems: Ergodic theory. \textit{Stochastics} \textbf{14} (1985), 273--317.

\bibitem{CoaylaTeranMohammedRuffino2007} E. A. Coayla-Teran, S.-E. A.
Mohammed and P. R. C. Ruffino. Hartman-Grobman theorems along hyperbolic
stationary trajectories. \textit{Discrete Contin. Dyn. Syst.} \textbf{17}
(2007), 281--292.

\bibitem{CoaylaTeranRuffino2004} E. A. Coayla-Teran and P. R. C. Ruffino.
Stochastic versions of Hartman-Grobman theorems. \textit{Stoch. Dyn.}
\textbf{4} (2004), 571--593.

\bibitem{Cong1996} N. D. Cong. Topological classification of linear
hyperbolic cocycles. \textit{J. Dynam. Differential Equations} \textbf{8}
(1996), 427--467.

\bibitem{DellacherieMeyer1978} C. Dellacherie and P.-A. Meyer. \textit{Probabilities and Potential}. North Holland, 1978.

\bibitem{DemetriusGundlach1999} L. Demetrius and V. M. Gundlach.
Evolutionary dynamics in random environments. In: H. Crauel and V. M.
Gundlach (eds.) \textit{Stochastic Dynamics}. Springer, 1999, pp. 371--394.

\bibitem{Evstigneev1974} I. V. Evstigneev. Positive matrix-valued cocycles
over dynamical systems. \textit{Uspekhi Matem. Nauk (Russ. Math. Surveys)}
\textbf{29} No. 5 (1974), 219--220. (In Russian.)

\bibitem{EvstigneevHensSchenkHoppe2006} I. V. Evstigneev, T. Hens and K. R.
Schenk-Hopp\'{e}. Evolutionary stable stock markets. \textit{Econom. Theory}
\textbf{27} (2006), 449--468.

\bibitem{EvstigneevHensSchenkHoppe2009} I. V. Evstigneev, T. Hens and K. R.
Schenk-Hopp\'{e}. Evolutionary finance. In: T. Hens and K. R. Schenk-Hopp\'{e} (eds.) \textit{Handbook of Financial Markets: Dynamics and Evolution}. North-Holland, 2009, pp. 507--566.

\bibitem{EvstigneevPirogov2007} I. V. Evstigneev and S. A. Pirogov. A
stochastic contraction principle. \textit{Random Oper. Stoch. Equ.} \textbf{15} (2007), 155--162.

\bibitem{EvstigneevPirogov2010} I. V. Evstigneev and S. A. Pirogov.
Stochastic nonlinear Perron-Frobenius theorem. \textit{Positivity} \textbf{14} (2010), 43--57.

\bibitem{EvstigneevSchenkHoppe2008} I. V. Evstigneev and K. R. Schenk-Hopp\'{e}. Stochastic equilibria in von Neumann-Gale dynamical systems. \textit{Trans. Amer. Math. Soc.} \textbf{360} (2008), 3345--3364.

\bibitem{FurstenbergKesten1960} H. Furstenberg and H. Kesten. Products of
random matrices. \textit{Ann. Math. Statist.} \textbf{31} (1960), 457--469.

\bibitem{Grobman1959} D. M. Grobman. Homeomorphisms of systems of
differential equations. \textit{Dokl. Akad. Nauk SSSR} \textbf{128} (1959),
880--881.

\bibitem{Hartman1960} P. Hartman. A lemma in the theory of structural
stability of differential equations. \textit{Proc. Amer. Math. Soc.} \textbf{11} (1960), 610--620.

\bibitem{Hartman1960a} P. Hartman. On local homeomorphisms of Euclidean
spaces. \textit{Bol. Soc. Mat. Mexicana} \textbf{5} (1960), 220--241.

\bibitem{KhaninKifer1996} K. Khanin and Yu. Kifer. Thermodynamic formalism
for random transformations and statistical mechanics. \textit{Amer. Math.
Soc. Transl. Ser. 2} \textbf{171} (1996), 107--140.

\bibitem{Kifer1996a} Yu. Kifer. Fractal dimensions and random
transformations. \textit{Trans. Amer. Math. Soc.} \textbf{348} (1996),
2003--2038.

\bibitem{Kifer1996b} Yu. Kifer. Perron-Frobenius theorem, large deviations,
and random perturbations in random environments. \textit{Math. Z.} \textbf{222} (1996), 677--698.

\bibitem{Kifer1998} Yu. Kifer. Limit theorems for random transformations and
processes in random environments. \textit{Trans. Amer. Math. Soc.} \textbf{350} (1998), 1481--1518.

\bibitem{Kifer2009} Yu. Kifer. Thermodynamic formalism for random
transformations revisited. \textit{Stoch. Dyn.} \textbf{8} (2008), 77--102.

\bibitem{KiferLyu2006} Yu. Kifer and P.-D. Liu. Random dynamics. In: B.
Hasselblatt and A. Katok (eds.) \textit{Handbook of Dynamical Systems}, Vol.
1B. Elsevier, 2006, pp. 379--499.

\bibitem{KolmogorovFomin1957} A. N. Kolmogorov and S. V. Fomin. \textit{Elements of the Theory of Functions and Functional Analysis}. Graylock, 1957.

\bibitem{LianLu2010} Z. Lian and K. Lu. Lyapunov exponents and invariant
manifolds for random dynamical systems in a Banach space. \textit{Mem.
Amer. Math. Soc.} \textbf{206} (2010), no. 967.

\bibitem{MacLeanThorpZiemba2011} L. C. MacLean, E. O. Thorp and W. T. Ziemba
(eds.). \textit{The Kelly Capital Growth Investment Criterion: Theory and
Practice}. World Scientific, 2011.

\bibitem{MohammedScheutzow1999} S.-E. A. Mohammed and M. K. R. Scheutzow. The
stable manifold theorem for stochastic differential equations. \textit{Ann.
Probab.} \textbf{27} (1999), 615--652.

\bibitem{MohammedZhangZhao2008} S.-E. A. Mohammed, T. Zhang and H. Zhao. The
stable manifold theorem for semilinear stochastic evolution equations and
stochastic partial differential equations. \textit{Mem. Amer. Math.
Soc.} \textbf{196} (2008), no. 917.

\bibitem{OchsOseledets1999} G. Ochs and V. I. Oseledets. Topological fixed
point theorems do not hold for random dynamical systems. \textit{J. Dynam.
Differential Equations} \textbf{11} (1999), 583--593.

\bibitem{Oseledets1968} V. I. Oseledets. A multiplicative ergodic theorem:
Lyapunov characteristic numbers for dynamical systems. \textit{Trans. Moscow
Math. Soc.} \textbf{19} (1968), 197--231.

\bibitem{PirogovSinai19751976} S. A. Pirogov and Ya. G. Sinai. Phase diagrams of classical lattice systems, I and II, \textit{Theor. Math. Phys.} \textbf{25} (1975), 1185--1192; III, \textit{Theor. Math. Phys.} \textbf{26} (1976), 39--49.

\bibitem{Ruelle1979} D. Ruelle. Ergodic theory of differentiable dynamical
systems. \textit{Publ. Math. Inst. Hautes {\'E}tudes Sci.} \textbf{50}
(1979), 27--58.

\bibitem{Ruelle1982} D. Ruelle. Characteristic exponents and invariant
manifolds in Hilbert space. \textit{Ann. of Math.} \textbf{115} (1982),
243--290.

\bibitem{Schmalfuss1998} B. Schmalfuss. A random fixed point theorem and the
random graph transformation. \textit{J. Math. Anal. Appl.}  \textbf{225} (1998), 91--113.

\bibitem{Sinai1972} Ya. G. Sinai. Gibbs measures in ergodic theory. \textit{Uspekhi Matem. Nauk (Russ. Math. Surveys)} \textbf{27} (1972), 21--69.

\bibitem{Wanner1994} T. Wanner. Linearization of random dynamical systems.
In: C. John, U. Kirchgraber and H. O. Walther (eds.) \textit{Dynamics Reported}, Vol. 4. Springer, 1994, pp. 203--269.
\end{thebibliography}
\end{document}